\documentclass[10pt,a4paper]{article}
\usepackage[utf8]{inputenc}
\usepackage{amsmath}
\usepackage{amsfonts}
\usepackage{amsmath}
\usepackage{amssymb}
\usepackage{amsthm}
\usepackage{mathrsfs}
\usepackage{enumerate}
\usepackage{mathptmx}
\usepackage{hyperref}
\usepackage{scalerel}
\usepackage{setspace}
\usepackage[scr=rsfso]{mathalfa}
\usepackage[Symbol]{upgreek}

\usepackage{fancyhdr}
\usepackage{etoolbox}
\usepackage[T1]{fontenc}

\DeclareSymbolFont{Symbols}{OMS}{zplm}{m}{n}
\DeclareMathSymbol{\infinity}{\mathord}{Symbols}{"31}
\DeclareMathAlphabet{\mathcal}{OMS}{cmsy}{m}{n}

\allowdisplaybreaks

\usepackage[a4paper, right=4.0cm, left=4.0cm, top=4.4cm, bottom=4.4cm]{geometry}

\usepackage{titlesec}
\titleformat{\section}{\large\bfseries}{\thesection}{1em}{}
\numberwithin{equation}{section}

\titleformat{\subsection}[runin]
{\normalfont\bfseries}{\thesubsection}{1em}{}

\newtheoremstyle{noenddot}
{\topsep}
{\topsep}
{\itshape}
{-4pt}
{\bfseries}
{}
{ }
{\thmname{#1}\thmnumber{ #2}\thmnote{ \normalfont(#3)}}

\theoremstyle{noenddot}
\newtheorem{theorem}{}[section]

\newcommand{\thm}{\textnormal{\textbf{Theorem.~}}}
\newcommand{\proposition}{\textnormal{\textbf{Proposition.~}}}
\newcommand{\corollary}{\textnormal{\textbf{Corollary.~}}}
\newcommand{\lemma}{\textnormal{\textbf{Lemma.~}}}

\newcommand{\newalpha}{\scaleobj{0.9}{\alpha}}
\newcommand{\newbeta}{\scaleobj{0.9}{\beta}}
\newcommand{\newmu}{\scaleobj{0.9}{\mu}}

\newcommand{\newtheta}{\scaleobj{0.92}{\theta}}
\newcommand{\newlambda}{\scaleobj{0.9}{\lambda}}
\newcommand{\newsigma}{\scaleobj{0.9}{\sigma}}
\newcommand{\newomega}{\scaleobj{0.9}{\omega}}

\newcommand{\newgamma}{\scaleobj{0.9}{\gamma}}
\newcommand{\newtau}{\scaleobj{0.9}{\tau}}
\newcommand{\neweta}{\scaleobj{0.9}{\eta}}
\newcommand{\newrho}{\scaleobj{0.9}{\rho}}

\newcommand{\newdelta}{\scaleobj{0.9}{\delta}}

\newcommand{\newxi}{\scaleobj{0.9}{\xi}}
\newcommand{\newzeta}{\scaleobj{0.9}{\zeta}}

\newcommand{\newupphi}{\scaleobj{0.95}{\upphi}}
\newcommand{\newuppsi}{\scaleobj{0.9}{\uppsi}}
\newcommand{\newuppi}{\scaleobj{0.9}{\uppi}}

\newcommand{\varin}{~\!\scaleobj{0.9}{\in}~\!}
\newcommand{\newpartial}{\scaleobj{0.95}{\partial}}
\newcommand{\varni}{~\!\scaleobj{0.9}{\ni}~\!}

\newcommand{\1}{1\!\!\!1}

\title{\vspace{-5ex}\large\textbf{EXPONENTIAL $\Phi$-MIXING IMPLIES EXPONENTIAL $\Psi$-MIXING FOR MARKOV FIELDS ON BOUNDED DEGREE GRAPHS}}
\author{\normalsize\textsc{Elias Zimmermann}}
\date{\vspace{-5ex}}

\begin{document}
\begin{spacing}{1.15}
    \maketitle

    \renewcommand{\abstractname}{\vspace{-\baselineskip}}

    \begin{abstract}
        \textsc{Abstract}. We show that for non-degenerate $k$-Markovian random fields with finite state space over a bounded degree graph with exponential growth rate $\newtheta$ uniform $\newupphi$-mixing with exponential decay rate $\newlambda > 3\newtheta$ implies uniform $\newuppsi$-mixing with exponential decay rate $(\newlambda {-} 3\newtheta)/9$. As an application we obtain exponential $\newuppsi$-mixing for Gibbs fields on regular trees arising from finite range potentials such as the Ising model at low inverse temperature or the Potts model with sufficiently many spin states. 
    \end{abstract}

    \section{Introduction}

    The aim of this paper is to study quantitative mixing properties of random fields over graphs. In this context mixing describes the phenomenon that the stochastic dependence of events taking place in sufficiently distant regions of the graph becomes arbitrarily small. A classical approach to mixing properties is the study of so called mixing coefficients. To be more precise let us consider a countably infinite, connected and simple graph $\mathscr{G} = (S,B)$ with sites $S$ and bonds $B$. Considering a probability space $(\Omega,\mathscr{F}\!\!,\newmu)$ together with a finite set of states $E$ let $(\newsigma_{s})_{s \!\!\varin\!\! S}$ be a family of random variables $\newsigma_{s}\colon \Omega \to E$. Let $\newdelta$ be a given measure of stochastic independence of subsets of $\Omega$. Then the mixing coefficient corresponding to $\newdelta$ is given by the map $\newuppi$ assigning to any two subsets $U,V \subseteq S$ the supremum
        \[\newuppi(U,V) := \sup\big\{\newdelta(A,B)\colon A \varin \mathscr{F}_{~\!\!U}, B \varin \mathscr{F}_{~\!\!V}\big\},\]
    where $\mathscr{F}_{~\!\!U}$ and $\mathscr{F}_{~\!\!V}$ denote the $\newsigma$-algebras generated by the random variables $(\newsigma_{s})_{s \!\!\varin\!\! U}$ and $(\newsigma_{s})_{s \!\!\varin\!\! V}$ respectively. For a detailed exposition of the topic of mixing coefficients in the context of random fields we refer to \cite{Dou94} and \cite{Bra07}.
    
    In many contexts one is interested in quantitative upper bounds for the values of a mixing coefficient $\newuppi$ in terms of the size and the distance of the sets $U$ and $V$. Below we will introduce a notion of exponential $\newuppi$-mixing, which implies in particular that for sufficiently distant finite sets $U,V \subseteq S$ we have 
        \[\newuppi(U,V) \leq C~\!e^{-\newlambda d(U,V)}\]
    for some constants $C$ and $\newlambda$, where $d$ denotes the path metric on the graph $\mathscr{G}$. In what follows we will refer to $\newlambda$ as the \textit{exponential decay rate} of $\newuppi$-mixing. In this paper we will be specifically interested in the mixing coefficients $\newupphi$ and $\newuppsi$, which correspond to the measures of stochastic independence $\newdelta_{\newupphi}$ and $\newdelta_{\newuppsi}$ defined by 
        \[\newdelta_{\newupphi}(A,B) = |\newmu(A|B) - \newmu(A)|\]
    and 
        \[\newdelta_{\newuppsi}(A,B) = \left|\frac{\newmu(A \cap B)}{\newmu(A)\newmu(B)} - 1 \right|\]
    for sets $A,B \varin \mathscr{F}$ of positive measure respectively. While it is straightforward to see that exponential $\newupphi$-mixing implies exponential $\newuppsi$-mixing with the same rate of decay, exponential $\newuppsi$-mixing appears to be much stronger a priori. Both mixing conditions have been studied in various contexts and are intimately linked to laws of large numbers and central limit theorems, see \cite{Bra05} and \cite{Bra07} and the references given therein. In \cite{PZ25+} it is shown that exponential $\newuppsi$-mixing implies a spherical equipartition property for automorphism invariant random fields on regular trees.
    
    The main result of this paper shows that for a non-degenerate $k$-Markovian random field $(\newsigma_{s})_{s \!\varin\! S}$ on a bounded degree graph $\mathscr{G}$ of exponential growth rate $\newtheta$ the following holds true: If $(\newsigma_{s})_{s \!\varin\! S}$ is exponential $\newupphi$-mixing with decay rate $\newlambda > 3\newtheta$, then $(\newsigma_{s})_{s \!\varin\! S}$ is exponential $\newuppsi$-mixing with decay rate $(\newlambda{-}3\newtheta)/9$. This generalizes and refines a result of Alexander for random fields over the $d$-dimensional square lattice $\mathscr{L}_{d}$, see \cite{Ale98}. Note that in the latter paper $\newupphi$-mixing is called weak mixing and $\newuppsi$-mixing is called ratio weak mixing. Since $\mathscr{L}_{d}$ has exponential growth rate $\newtheta = 0$ the result is true here for every decay rate $\newlambda > 0$. Based on several results regarding exponential $\newupphi$-mixing of Gibbs fields over square lattices, Alexander used this result to obtain exponential $\newuppsi$-mixing for various models throughout the uniqueness region. While for general graphs similarly fine results on $\newupphi$-mixing are not at hand, it is well known that Gibbs fields arising from finite range potentials with sufficiently small Dobrushin constant satisfy exponential $\newupphi$-mixing with arbitrarily high decay rate. Thus one can use our main result to determine examples of exponentially $\newuppsi$-mixing Gibbs fields beyond square lattices. We illustrate this in the case of the $d$-regular tree $\mathscr{T}_{d}$, where we obtain exponential $\newuppsi$-mixing for the Ising model  at low inverse temperature as well as for Potts models with many spin states. We also give explicit estimates for the rate of decay in terms of the parameters of the models. 
    
        \medskip
        
    The paper is organized as follows: In Section \ref{preliminaries} we introduce the relevant objects, concepts and notational conventions. Section \ref{main} is devoted to the proof of our main result, Theorem \ref{Thm1}. In Section \ref{applications} we illustrate applications of our main result in the case of Gibbs fields. 

    \medskip

    \noindent \textbf{Acknowledgements.} The paper is part of the author's PhD project. The author wants to thank his advisor Felix Pogorzelski for several helpful discussions and comments in various stages of the work as well as constant support during the writing process. He is also indebted to Amos Nevo for many inspiring conversations. The author gratefully acknowledges support of the German-Israeli foundation for Science and Development (GIF) through grant I-1485/304.6-2019 and the German Academic Scholarship Foundation through a doctoral scholarship. 

    \section{Preliminaries} \label{preliminaries}

    Let $(\Omega_{1},\mathscr{A}_{1})$ and $(\Omega_{2},\mathscr{A}_{2})$ be standard Borel spaces. By a \textit{probability kernel} from $\mathscr{A}_{1}$ to $\mathscr{A}_{2}$ we mean a map $\newrho\colon \Omega \times \mathscr{A}_{2} \to [0,1]$ such that for any $\newomega \varin \Omega$ the map $\newrho(\newomega,\cdot)$ is a probability measure, which will be denoted by $\newrho^{\newomega}$ in what follows, and for any $A \varin \mathscr{A}_{2}$ the map $\newrho(\cdot,A)$ is $\mathscr{A}_{1}$-measurable. Let $\newmu_{1}$ be a probability measure on $(\Omega_{1},\mathscr{A}_{1})$ and $\newmu_{2}$ be a probability measure on $(\Omega_{2},\mathscr{A}_{2})$ and let for $i = 1,2$ the map $\newsigma_{i}\colon \Omega_{1} \times \Omega_{2} \to \Omega_{i}$ be the projection to the coordinate $i$. By a \textit{coupling} of $\newmu_{1}$ and $\newmu_{2}$ we mean a probability measure $\neweta$ on $(\Omega_{1} \times \Omega_{2},\mathscr{A}_{1} \times \mathscr{A}_{2})$ such that ${\newsigma_{1}}_{*}\neweta = \newmu_{1}$ and ${\newsigma_{2}}_{*}\neweta = \newmu_{2}$.

Throughout this paper $\mathscr{G}$ will always be a countably infinite and simple (undirected) graph, i.\@ e. $\mathscr{G}$ is given by a pair $(S,B)$ consisting of a countably infinite set $S$ of \textit{sites} and an antireflexive (symmetric) set $B \subseteq S \times S$ of \textit{bonds}. In what follows we will often write $s \sim t$ instead of $(s,t) \varin B$. Given a site $s \varin S$ the \textit{degree} $\mathrm{deg}(s)$ is given by the number of sites $t \varin S$ with $s \sim t$. We shall say that $\mathscr{G}$ is of \textit{bounded degree} if there is a constant $d > 0$ such that $\mathrm{deg}(s) \leq d$ for every $s \varin S$. A sequence $s_{1},\ldots,s_{m}$ of sites $s_{i} \varin S$ is called a \textit{path} if it satisfies $s_{i} \sim s_{i+1}$ for $i = 1,\ldots,m{-}1$. The \textit{length} of the path is given by $m{-}1$. We say that $\mathscr{G}$ is \textit{connected} if for any sites $s,t \varin S$ there is a path $r_{1},\ldots,r_{m}$ with $r_{1} = s$ and $r_{m} = t$. 

In what follows we shall always assume that $\mathscr{G}$ is connected and of bounded degree. In this case we obtain a metric $d$ on $S$ by assigning to every pair of sites $s,t \varin S$ the length $d(s,t)$ of a minimal path connecting $s$ and $t$. The metric $d$ is called the \textit{path metric} of $\mathscr{G}$. As usual $d$ extends to a distance function on subsets of $S$ given by 
    \[d(U,V) := \min\{d(s,t)\colon s \varin U, t \varin V\}\]
for $U,V \subseteq S$. Given a site $s \varin S$ and a set $\Lambda \subseteq S$ we shall write $d(s,\Lambda)$ for the distance of the sets $\{s\}$ and $\Lambda$. Fixing $k \in \mathbb{N}_{0}$ we define the \textit{$k$-boundary} of a set $\Lambda \subseteq S$ with respect to $d$  by
    \[\newpartial_{\!k}\Lambda := \{s \varin \Lambda^{c}\colon d(s,\Lambda) \leq k\}.\]
Moreover, fixing $s \varin S$ and $n \in \mathbb{N}_{0}$, we shall denote by 
    \[B_{n}(s)  := \big\{t \varin S\colon d(s,t) \leq n\big\}\]
the \textit{ball} of radius $n$ around $s$ with respect to $d$. It is not difficult to see that for a bounded degree graph $\mathscr{G}$ there is a $\newtheta \varin [0,\infinity)$ such that $|B_{n}(s)| \leq e^{\newtheta n}$ for every $s \varin S$ and $n \varin \mathbb{N}$. We shall refer to the infimum of all such $\newtheta$ as the \textit{exponential growth rate} of $\mathscr{G}$.
    
Let $E = \{1,\ldots,k\}$ be a finite set of states. Let us denote by $\Omega := E^{S}$ the respective set of configurations on $S$. Then for every site $s \varin S$ we may consider the projection $\newsigma_{s}\colon \Omega \to E$ given by $\newsigma_{s}(\newomega) = \newomega_{s}$ for $\newomega \varin \Omega$. More generally, we shall denote by $\newsigma_{\Lambda}$ the projection mapping a configuration $\newomega$ to the restricted configuration $\newomega_{\Lambda}$ for every set $\Lambda \subseteq S$. The topology generated by the family of projections $\{\newsigma_{s}\colon s \varin S\}$ turns $\Omega$ into a compact Polish space. Accordingly, equipping $\Omega$ with the respective Borel $\newsigma$-algebra $\mathscr{F}$, we obtain a standard Borel space $(\Omega,\mathscr{F})$. Let us introduce some notation: For any subset $\Lambda \subseteq S$ we shall write $\mathscr{F}_{\Lambda}$ for the $\newsigma$-algebra generated by the projections $\{\newsigma_{s}\colon s \varin \Lambda\}$. Given two sets $\Delta,\Lambda \subseteq S$ we shall write $\Delta \sqsubseteq \Lambda$ to indicate that $\Delta$ is a finite subset of $\Lambda$. Moreover, given a subset $\Lambda \sqsubseteq S$ and a configuration $u \varin E^{\Lambda}$ we shall denote by $[u]
$ the cylinder set consisting of all $\newomega \varin \Omega$ with $u = \newomega_{\Lambda}$. Considering a probability measure $\newmu$ on $(\Omega,\mathscr{F})$ the family of projections $\{\newsigma_{s}\colon s \varin S\}$ constitutes a random field. On the other hand, it is easy to see that any random field over $S$ taking values in $E$ is equivalent in the distributional sense to a random field of the above form. Thus, there is no restriction of generality in confining to random fields of the latter kind. 

Since $(\Omega,\mathscr{F})$ is a standard Borel space, any probability measure $\newmu$ on $(\Omega,\mathscr{F})$ admits \textit{regular conditional probabilities}, i.\@ e. there is an a.\@ s. unique family $(\newmu_{\Lambda})_{\Lambda \sqsubseteq S}$ of probability kernels $\newmu_{\Lambda}$ from $\mathscr{F}_{\Lambda^{c}}$ to $\mathscr{F}_{\Lambda}$ such that for every $\Lambda \sqsubseteq S$ and every $A \varin \mathscr{F}_{\Lambda}$ we have
    \[\newmu^{\newomega}_{\Lambda}(A) = \mathbb{E}[\1_{A}|\mathscr{F}_{\Lambda^{c}}](\newomega)\]
for $\newmu$-almost all $\newomega \varin \Omega$. We shall say that the probability measure $\newmu$ is  \textit{non-degenerate} if for every $\Lambda \sqsubseteq S$ and every $A \varin \mathscr{F}_{\Lambda}$ with $\newmu(A) > 0$ one has $\newmu_{\Lambda}^{\omega}(A) > 0$ for $\newmu$-almost all $\newomega \varin \Omega$. Given $k \varin \mathbb{N}_{0}$ a probability measure $\newmu$ (as well as the respective random field) is called \textit{$k$-Markovian} if $\newomega_{\newpartial_{k}\Lambda} = \newtau_{\newpartial_{k}\Lambda}$ implies $\newmu_{\Lambda}^{\newomega} = \newmu_{\Lambda}^{\tau}$ for $\newmu$-almost all $\newomega,\newtau \varin \Omega$. 

\section{The main result} \label{main}

This section is devoted to the proof of the main result. Throughout this section $\mathscr{G} = (S,B)$ will be a countable connected graph of bounded degree. We shall denote by $\newtheta$ the exponential growth rate of $\mathscr{G}$. Let $(\newsigma_{s})_{s \!\varin\! S}$ be a random field over $\mathscr{G}$ taking values in a finite set $E$. We assume that $(\newsigma_{s})_{s \!\varin\! S}$ is given in canonical form. Let $(\Omega,\mathscr{F}\!,\newmu)$ be the corresponding space of configurations. We will also assume that $\newmu$ is non-degenerate. For any two sets $U,V \subseteq S$ together with a parameter $\newlambda > 0$ we shall set
    \[\Sigma(\newlambda,U,V) := \sum_{t \!\varin\! U}\sum_{s \!\varin\! V}e^{-\newlambda d(t,s)}.\]
Given a mixing coefficient $\newuppi$ we will call the measure $\newmu$ (and the respective random field) \textit{$\newuppi$-mixing with exponential decay rate $\newlambda$} if there is a constant $C \geq 1$ such that for all $U,V \subseteq S$ with $C~\Sigma(\newlambda,U,V) < 1$ we have
    \begin{align} \label{eq11}
    \newuppi(U,V)  \leq C ~ \Sigma(\newlambda,U,V).
    \end{align}
Note that the condition $C~\Sigma(\newlambda,U,V) < 1$ implies in particular that the sets $U$ and $V$ are disjoint. If $\newuppi$ equals the coefficient $\newupphi$, we have $\newuppi(U,V) \leq 1$ by definition, so in this case exponential mixing implies the validity of (\ref{eq11}) for all disjoint sets $U,V \subseteq S$. Note that if $U$ and $V$ are finite sets $(\ref{eq11})$ yields the upper bound
    \[\newuppi(U,V) \leq C~\!|U|~\!|V|~\!e^{-\newlambda d(U,V)},\]
which imposes a tradeoff between the size and the distance of the sets $U, V$. More precisely it states that the stochastic dependence of events taking place in $U$ and $V$, as measured by $\newuppi$, grows at most linearly in the sizes and decreases at least exponentially in the distance of $U$ and $V$. 

We will start by proving a series of technical statements. The first one is a characterization of exponential $\newupphi$-mixing, which is well known for random fields over square lattices. Note that in this case the exponential growth rate is $0$, so the statement holds for every decay rate $\newlambda > 0$.

\begin{theorem} \label{prop1} \proposition Fix $\newlambda > 3\newtheta$. Then $\newmu$ is $\newupphi$-mixing with exponential decay rate $\newlambda$ if and only if there is a constant $C_{1} \geq 1$ such that for all sets $\Delta \subseteq \Lambda \sqsubseteq S$ the inequality
    \begin{align}
        \big\|{\newsigma_{\Delta}}_{*}\newmu_{\Lambda}^{\newomega} - {\newsigma_{\Delta}}_{*}\newmu_{\Lambda}^{\newtau}\big\|_{\textnormal{Var}} \leq C_{1}~\Sigma(\newlambda,\Delta,\Lambda^{c})
    \end{align}
holds for $\newmu$-almost all $\newomega,\newtau \varin \Omega$.
\end{theorem}

\textit{Proof:} To prove the only-if-direction let $C > 0$ be the constant witnessing the $\newlambda$-exponential $\newupphi$-mixing of $\newmu$. Fix $\Delta \subseteq \Lambda \sqsubseteq S$ and $A \varin \mathscr{F}_{\Delta}$ with $\newmu(A) > 0$. Consider the sequence $(\Lambda_{n})_{n=1}^{\infinity}$ consisting of the sets $\Lambda_{n} := \Lambda \cup B_{n}^{c}$. Note that for $\newmu$-almost all $\newomega \in \Omega$ and all $n \in \mathbb{N}$ the cylinder set $[\newomega_{\Lambda_{n}^{c}}]$ has positive measure. (In fact, since $\Omega$ is a Polish space, the support of $\newmu$ whose elements obviously satisfy the above property, has full measure.) Thus, observing that $\Lambda^{c} = \bigcup_{n=1}^{\infinity} \Lambda_{n}^{c}$, we may use Levy's $0$-$1$-law to obtain that
    \begin{align*}
        |\newmu_{\Lambda}^{\newomega}(A) - \newmu(A)|  &= \big|\mathbb{E}[\1_{A}|\mathscr{F}_{\Lambda^{c}}](\newomega) - \newmu(A)\big| = \lim_{n \to \infinity}\big|\mathbb{E}[\1_{A}|\mathscr{F}_{\Lambda_{n}^{c}}](\newomega) - \newmu(A)\big| \\
    &= \lim_{n \to \infinity}|\newmu(A|\newsigma_{\Lambda_{n}^{c}} = \newomega_{\Lambda_{n}^{c}}) - \newmu(A)|
     \leq \limsup_{n \to \infinity}~ C~\Sigma(\newlambda,\Delta,\Lambda_{n}^{c}) \\
&\leq C~\Sigma(\newlambda,\Delta,\Lambda^{c})
    \end{align*}
for $\newmu$-almost $\newomega \in \Omega$. Since $\mathscr{F}_{\Delta}$ contains only finitely many sets we may conclude that for all $\newomega,\newtau$ in a set of full measure we have
    \begin{align*}
        |\newmu_{\Lambda}^{\newomega}(A) - \newmu^{\newtau}_{\Lambda}(A)| \leq |\newmu_{\Lambda}^{\newomega}(A) - \newmu(A)| + |\newmu(A) - \newmu_{\Lambda}^{\newtau}(A)| \le 2~C~\Sigma(\newlambda,\Delta,\Lambda^{c}).
    \end{align*}
This shows the claim.

To verify the other direction fix again $\Delta \subseteq \Lambda \sqsubseteq T$. Then for all $A \in \mathscr{F}_{\Delta}$ and $B \in \mathscr{F}_{\Lambda^{c}}$ with $\newmu(A)\newmu(B) > 0$ we have
    \begin{align} \label{eq4}
        \newmu(A \cap B) = \int_{B}\1_{A}(\newomega)~d\newmu(\newomega) = \int_{B} \mathbb{E}[\1_{A}|\mathscr{F}_{\Lambda^{c}}](\newomega) ~d\newmu(\newomega) = \int_{B}\newmu_{\Lambda}^{\newomega}(A)~d\newmu(\newomega),
    \end{align}
which implies
    \begin{align*}
        |\newmu(A|B) - \newmu(A)| &= \bigg|\frac{1}{\newmu(B)}\int_{B}\newmu_{\Lambda}^{\newomega}(A)~d\newmu(\newomega) - \int_{\Omega}\newmu_{\Lambda}^{\newtau}(A)~d\newmu(\newtau)\bigg| \\
&\leq \frac{1}{\newmu(B)}\int_{B}\int_{\Omega}|\newmu_{\Lambda}^{\newomega}(A) - \newmu_{\Lambda}^{\newtau}(A)| ~d\newmu(\newtau)d\newmu(\newomega) \\ &\le C_{1}~\Sigma(\newlambda,\Delta,\Lambda^{c}). \vphantom{\int}
    \end{align*}
Given arbitrary disjoint sets $U,V \subseteq S$ we consider the sequences $(\Delta_{n})_{n=1}^{\infinity}$ and $(\Lambda_{n})_{n=1}^{\infinity}$ defined by $\Delta_{n} := U \cap B_{n}$ and $\Lambda_{n} := V^{c} \cap B_{2n}$. Then we have $\Delta_{n} \subseteq \Lambda_{n} \sqsubseteq T$ for all $n \in \mathbb{N}$ and obtain 
    \begin{align*}
        \Sigma(\newlambda,\Delta_{n},\Lambda_{n}^{c}) &= \sum_{s\!\!\varin \!\!\Delta_{n}}\sum_{t\!\!\varin \!\!\Lambda_{n}^{c}} e^{-\newlambda d(s,t)} \leq \sum_{s\!\!\varin \!\!\Delta_{n}}\sum_{t\!\!\varin \!\!V} e^{-\newlambda d(s,t)} + \sum_{s\!\!\varin \!\!\Delta_{n}}\sum_{t\!\!\varin \!\! B_{2n}^{c}} e^{-\newlambda d(s,t)} \\
&\le \sum_{s \!\!\varin \!\! U}\sum_{t \!\!\varin \!\! V} e^{-\newlambda d(s,t)} + \sum_{s \!\!\varin \!\! B_{n}}\sum_{t \!\!\varin \!\! B_{2n}^{c}} e^{-\newlambda d(s,t)}.
    \end{align*}
The first term in the above expression equals $\Sigma(\newlambda,U,V)$. Setting $\newbeta := \exp(\newlambda - \newtheta)$ we may estimate the second term from above by
    \begin{align*}
        \sum_{s \!\!\varin \!\! B_{n}}\sum_{t \!\!\varin \!\! B_{2n}^{c}} e^{-\newlambda(|t|-|s|)}&\leq |B_{n}|\sum_{t \!\!\varin \!\! B_{2n}^{c}} e^{-\newlambda(|t| - n)} = |B_{n}| \sum_{m= 2n+1}^{\infinity} |S_{m}| e^{-\newlambda(m-n)} \\
&= |B_{n}|\sum_{m = n+1}^{\infinity}|S_{m+n}|e^{-\newlambda m} \leq e^{\newtheta n}\sum_{m=n}^{\infinity} e^{\newtheta (m+n)} e^{-\newlambda m} \\
&= e^{2\newtheta n}\sum_{m=n}^{\infinity}e^{m(\newtheta - \newlambda)} = \frac{\newbeta e^{n(3 \newtheta - \newlambda)}}{\newbeta {-} 1} .
    \end{align*}
Since the last term goes to zero for $n \to \infinity$ by the choice of $\newlambda$, we may conclude  that
    \begin{align} \label{eq2}
        \limsup_{n \to \infinity}~\Sigma(\newlambda, \Delta_{n},\Lambda_{n}^{c}) \leq \Sigma(\newlambda,U,V)
    \end{align}
Now fix sets $A \varin \mathscr{F}_{U}$ and $B \varin \mathscr{F}_{V}$ with $\newmu(A)\newmu(B) > 0$. Then $B \varin \mathscr{F}_{\Lambda_{n}^{c}}$ for all $n \in \mathbb{N}$. By a well known approximation lemma there exists a sequence $(A_{n})_{n=1}^{\infinity}$ of sets $A_{n} \in \mathscr{F}_{\Delta_{n}}$ with $\newmu(A_{n} \triangle A) \to 0$, which implies $\newmu(A_{n}|B) \to \newmu(A|B)$ as well as $\newmu(A_{n}) \to \newmu(A)$. Together with (\ref{eq2}) this gives
    \begin{align*}
        |\newmu(A|B) - \newmu(A)| &= \lim_{n \to \infinity}|\newmu(A_{n}|B) - \newmu(A_{n})| \\ &\leq \limsup_{n \to \infinity}~C_{1}~\!\Sigma(\newlambda,\Delta_{n},\Lambda_{n}^{c}) \\ &\leq C_{1}~\!\Sigma(\newlambda,U,V).
    \end{align*}
This shows that $\newmu$ is $\newupphi$-mixing with exponential decay rate $\newlambda$. \hfill $\Box$
\begin{theorem} \label{prop2} \proposition Fix $\newlambda > 3\newtheta$. Assume that there is a constant $C_{2} \geq 1$ such that for all sets $\Delta \subseteq \Lambda \sqsubseteq S$ with $C_{2}~\Sigma(\newlambda,\Delta,\Lambda^{c}) < 1$ and all $u_{\Delta} \varin E^{\Delta}$ satisfying $\newmu(\newsigma_{\Delta} = u_{\Delta}) > 0$ we have
      \begin{align} \label{eq3}
	\big|\newmu_{\Lambda}^{\newomega}(\newsigma_{\Delta} = u_{\Delta}) - \newmu_{\Lambda}^{\newtau}(\newsigma_{\Delta} = u_{\Delta})\big| \leq C_{2}~\Sigma(\newlambda,\Delta,\Lambda^{c})\newmu_{\Lambda}^{\newtau}(\newsigma_{\Delta} = u_{\Delta})
	\end{align}
for $\newmu$-almost all $\newomega,\newtau \varin \Omega$. Then $\newmu$ is $\newuppsi$-mixing with exponential decay rate $\lambda$.
\end{theorem}

\textit{Proof:} Fix $C_{2} \geq 1$ and $\Delta \subseteq \Lambda \sqsubseteq S$ as above. Consider a configuration $u_{\Delta} \varin E^{\Delta}$ with $\newmu(\newsigma_{\Delta} = u_{\Delta}) > 0$ and let $\Omega_{0} \subseteq \Omega$ be a set of full measure such that (\ref{eq3}) holds for all $\newomega, \newtau \in \Omega_{0}$. Fix $B \in \mathscr{F}_{\Lambda^{c}}$ with $\newmu(B) > 0$. Then by (\ref{eq4})we have
    \begin{align*}
        \big|\newmu(\newsigma_{\Delta} = u_{\Delta}|B) - \newmu(\newsigma_{\Delta} = &u_{\Delta})\big| = \bigg|\frac{1}{\newmu(B)}\int_{B}\newmu_{\Lambda}^{\newomega}(\newsigma_{\Delta} = u_{\Delta})~d\newmu(\newomega) - \int_{\Omega} \newmu_{\Lambda}^{\newtau}(\newsigma_{\Delta} = u_{\Delta})~d\newmu(\newtau)\bigg| \\
        &\leq \frac{1}{\newmu(B)}\int_{B \cap \Omega_{0}}\int_{\Omega_{0}}\big| \newmu_{\Lambda}^{\newomega}(\newsigma_{\Delta} = u_{\Delta}) - \newmu_{\Lambda}^{\newtau}(\newsigma_{\Delta} = u_{\Delta})\big|~d\newmu(\newtau)d\newmu(\newomega) \\
        &\leq C_{2}~\Sigma(\newlambda,\Delta,\Lambda^{c})\frac{\newmu(B \cap \Omega_{0})}{\newmu(B)} \int_{\Omega}\newmu_{\Lambda}^{\newtau}(\newsigma_{\Delta} = u_{\Delta})~d\newmu(\newtau) \\
        &= C_{2}~\Sigma(\newlambda,\Delta,\Lambda^{c}) \newmu(\newsigma_{\Delta} = u_{\Delta}). \vphantom{\int}
    \end{align*}
Now for a general $A \in \mathscr{F}_{\Delta}$ with $\newmu(A) > 0$ we can find configurations $u_{\Delta}^{1},\ldots,u_{\Delta}^{m} \varin E^{\Delta}$ with $\newmu(\newsigma_{\Delta} = u_{\Delta}^{i}) > 0$ for $i = 1,\ldots,m$ such that $A = [u_{\Delta}^{1}] \cup \ldots \cup [u_{\Delta}^{m}]$ up to a $\newmu$-null set. This implies
    \begin{align*}
        |\newmu(A|B) - \newmu(A)| &\leq \sum_{i=1}^{m} \big|\newmu(\newsigma_{\Delta} = u_{\Delta}^{i}|B) - \newmu(\newsigma_{\Delta} = u_{\Delta}^{i})\big| \\ &\leq C_{2}~\Sigma(\newlambda,\Delta,\Lambda^{c})\sum_{i=1}^{m}\newmu(\newsigma_{\Delta} = u_{\Delta}^{i}) \\ &= C_{2}~\Sigma(\newlambda,\Delta,\Lambda^{c})\newmu(A). \vphantom{\sum}
    \end{align*}
Finally, for $U, V \subseteq S$ arbitrary disjoint sets consider the sequences $(\Delta_{n})_{n=1}^{\infinity}$ and $(\Lambda_{n})_{n=1}^{\infinity}$ given by $\Delta_{n} := U \cap B_{n}$ and $\Lambda_{n} := V^{c} \cap B_{2n}$ for $n \in \mathbb{N}$. Then we have $\Delta_{n} \subseteq \Lambda_{n} \sqsubseteq S$ for all $n \in \mathbb{N}$. Fix $A \in \mathscr{F}_{~\!\!U}$ and $B \varin \mathscr{F}_{~\!\!V}$ with $\newmu(A)\newmu(B) > 0$. Then $B \varin \mathscr{F}_{\Lambda_{n}^{c}}$ for all $n \in \mathbb{N}$. Furthermore, as in the proof of Proposition \ref{prop1}, there is a sequence $(A_{n})_{n=1}^{\infinity}$ of sets $A_{n} \in \mathscr{F}_{\Delta_{n}}$ such that $\newmu(A_{n}|B) \to \newmu(A|B)$ and $\newmu(A_{n}) \to \newmu(A)$. As in the proof of Proposition \ref{prop1} the choice of $\newlambda$ gives
    \[\limsup_{n \to \infinity} \Sigma(\newlambda,\Delta_{n},\Lambda_{n}^{c}) \leq \Sigma(\newlambda,U,V),\]
which implies
    \begin{align*}
        |\newmu(A|B) - \newmu(A)| &= \lim_{n \to \infinity}|\newmu(A_{n}|B) - \newmu(A_{n})| \\ &\leq \limsup_{n \to \infinity}~ C_{2}~\Sigma(\newlambda,\Delta_{n},\Lambda_{n}^{c})\newmu(A_{n}) \\ &\leq C_{2}~\Sigma(\newlambda,\Delta,\Lambda^{c}) \newmu(A).
    \end{align*}
This shows that $\newmu$ is $\newuppsi$-mixing with exponential decay rate $\newlambda$. \hfill $\Box$

\begin{theorem} \label{lemma1}
    \lemma Assume that $\newmu$ is $\newupphi$-mixing with exponential decay rate $\newlambda > 0$. Then for some constant $C_{3} > 0$ the following holds. Fix $\Delta \subseteq \Lambda \sqsubseteq S$ and $\Upsilon \subseteq \Lambda\setminus \Delta$. Then for all $A \varin \mathscr{F}_{~\!\Upsilon}$ and $B \varin \mathscr{F}_{\Delta}$ with $\newmu(A)\newmu(B) > 0$ we have
        \[\big|\newmu_{\Lambda}^{\newomega}(A|B) - \newmu_{\Lambda}^{\newomega}(A)\big| \leq C_{3}~\Sigma(\newlambda,\Upsilon,\Delta \cup \Lambda^{c})\]
    for $\newmu$-almost all $\newomega \varin \Omega$.
\end{theorem}

\textit{Proof:} Let $C > 0$ be the constant witnessing the exponential $\newupphi$-mixing of $\newmu$. Consider the sequence $(\Lambda_{n})_{n=1}^{\infinity}$ of subsets $\Lambda_{n} \subseteq S$ given by $\Lambda_{n} := \Lambda \cup B_{n}^{c}$. Then we have $\Lambda^{c} := \bigcup_{n=1}^{\infinity}\Lambda_{n}^{c}$. Fix $A \varin \mathscr{F}_{~\!\Upsilon}$ and $B \varin \mathscr{F}_{\Delta}$ with $\newmu(A)\newmu(B) > 0$. Recall that for $\newmu$-almost all $\newomega \varin \Omega$ the cylinder set $[\newomega_{\Lambda_{n}^{c}}]$ has positive measure for all $n \in \mathbb{N}$. Thus, we may apply as above Levy's $0$-$1$-law to obtain
    \begin{align*}
        \big|\newmu_{\Lambda}^{\newomega}(A \cap B) - \newmu(A)\newmu_{\Lambda}^{\newomega}(B)\big|
        &= \lim_{n \to \infinity}\big|\newmu(A\cap B|\newsigma_{\Lambda_{n}^{c}} = \newomega_{\Lambda_{n}^{c}}) - \newmu(A)\newmu(B|\newsigma_{\Lambda_{n}^{c}} = \newomega_{\Lambda_{n}^{c}})\big| \\
        &= \lim_{n \to \infinity} \frac{\big|\newmu(A \cap B \cap [\newomega_{\Lambda_{n}^{c}}]) - \newmu(A)\newmu(B \cap [\newomega_{\Lambda_{n}^{c}}])\big|}{\newmu(\newsigma_{\Lambda_{n}^{c}} = \newomega_{\Lambda_{n}^{c}})} \\
        &\leq \limsup_{n \to \infinity}~C~\Sigma(\newlambda,\Upsilon,\Delta \cup \Lambda_{n}^{c})~\frac{\newmu(B \cap [\newomega_{\Lambda_{n}^{c}}])}{\newmu(\newsigma_{\Lambda_{n}^{c}} = \newomega_{\Lambda_{n}^{c}})} \\
        &\leq C~\Sigma(\newlambda,\Upsilon,\Delta \cup \Lambda^{c})~\lim_{n \to \infinity} \newmu(B|\newsigma_{\Lambda_{n}^{c}} = \newomega_{\Lambda_{n}^{c}}) \vphantom{\lim_{~}}\\
        &= C~\Sigma(\newlambda,\Upsilon,\Delta \cup \Lambda^{c})~\newmu_{\Lambda}^{\newomega}(B)
    \end{align*}
and thus
    \[\big|\newmu_{\Lambda}^{\newomega}(A|B) - \newmu(A)\big| \leq C~\Sigma(\newlambda,\Upsilon,\Delta \cup \Lambda^{c})\]
for $\newmu$-almost all $\newomega \varin \Omega$. Furthermore we have
    \begin{align*}
        \big|\newmu_{\Lambda}^{\newomega}(A) - \newmu(A)\big| &= \lim_{n \to \infinity}\big|\newmu(A|\newsigma_{\Lambda_{n}^{c}} = \newomega_{\Lambda_{n}^{c}}) - \newmu(A)\big| \leq \limsup_{n \to \infinity} C~\Sigma(\newlambda,\Upsilon,\Lambda_{n}^{c}) \\
        &\leq C~\Sigma(\newlambda,\Upsilon,\Delta \cup \Lambda^{c})
    \end{align*}
$\newmu$-almost all surely. Both together implies
    \[\big|\newmu_{\Lambda}^{\newomega}(A|B) - \newmu^{\newomega}_{\Lambda}(A)\big| \leq \big|\newmu_{\Lambda}^{\newomega}(A|B) - \newmu(A)\big| + \big|\newmu_{\Lambda}^{\newomega}(A) - \newmu(A)\big| \leq 2~\!C~\!\Sigma(\newlambda,\Upsilon,\Delta \cup \Lambda^{c})\]
for $\newmu$-almost all $\newomega \varin \Omega$. This shows the claim. \hfill $\Box$ \\

We are now able to prove our main result. In the special case that $\mathscr{G}$ equals the square lattice $\mathscr{L}_{d}$ it follows from a result of Alexander, see \cite[Theorem 3.3]{Ale98}. The proof presented here is generalization of Alexander's proof. The assumption on $\newlambda$ becomes necessary in order to handle exponential growth of balls in $\mathscr{G}$. In fact, Alexander showed the statement under a slightly weaker condition on the random field than being $k$-Markovian. However, since all examples we are going to discuss are $k$-Markovian and the necessary modifications of the proof to cover the more general setting are straightforward, we will stick to this assumption. For the proof we introduce the notation
    \[\newdelta(C,\newlambda,U,V) := C~\Sigma(\newlambda,U,V)\]
for any $C, \newlambda > 0$ and $U,V \subseteq S$.

\begin{theorem} \label{Thm1}
    \thm Fix $\newlambda > 3\newtheta$ and $k \in \mathbb{N}_{0}$. If $\newmu$ is $k$-Markovian and $\newupphi$-mixing with exponential decay rate $\newlambda$, then $\newmu$ is $\newuppsi$-mixing with exponential decay rate $(\newlambda {-} 3\newtheta)/9$. 
\end{theorem}

    \textit{Proof:} Let $C_{1}$ and $C_{3}$ be the constants appearing in Proposition \ref{prop1} and Lemma \ref{lemma1} and set $\newlambda_{0} := (\newlambda{-}3\newtheta)/9$, $C_{4} := \max\{1,C_{1},C_{3}\}$, $\newalpha := \exp(\newtheta{-} \newlambda/3)$ and $\newbeta := \newalpha/(1{-}\newalpha)$. Furthermore, choose $C := \max\{12C_{4}\newbeta,\exp(3k\newlambda_{0})\}$.
    Since there are only countably many finite subsets of $S$ we may fix a set $\Omega_{0} \subseteq \Omega$ of full measure such that the statement of Proposition \ref{prop1} holds for all $\Delta \subseteq \Lambda \sqsubseteq S$ and $\newomega,\newtau \in \Omega_{0}$ and the statement of Lemma \ref{lemma1} holds for all $\Delta \subseteq \Lambda \sqsubseteq S$, $\Upsilon \subseteq \Lambda\setminus\Delta$ and $\newomega \varin \Omega_{0}$. Since $\newmu$ is non-degenerate we may furthermore choose $\Omega_{0}$ in such a way that in addition for all sets $\Delta \subseteq \Lambda \sqsubseteq S$, all $\newomega \in \Omega_{0}$ and any $u_{\Delta} \in E^{\Delta}$ with $\newmu(\newsigma_{\Delta} = u_{\Delta}) > 0$ we have $\newmu_{\Lambda}^{\newomega}(\newsigma_{\Delta} = u_{\Delta}) > 0$.

    By Proposition \ref{prop2} it suffices to show that for all $\Delta \subseteq \Lambda \sqsubseteq S$ with $\newdelta(C,\newlambda_{0},\Delta,\Lambda^{c}) < 1$ and all configurations $u^{~}_{\Delta} \varin E^{\Delta}$ with $\newmu(\newsigma_{\Delta} = u_{\Delta}) > 0$ we have
        \begin{align} \label{eq5}
            \big|\newmu_{\Lambda}^{\newomega}(\newsigma_{\Delta} = u_{\Delta}) - \newmu_{\Lambda}^{\newtau}(\newsigma_{\Delta} = u_{\Delta})\big| \leq \newdelta(C,\newlambda_{0},\Delta,\Lambda^{c})~\newmu_{\Lambda}^{\newtau}(\newsigma_{\Delta} = u_{\Delta})
        \end{align}
    for all $\newomega,\newtau \varin \Omega_{0}$. To show this fix $\newomega,\newtau \varin \Omega_{0}$ and $\Delta,\Lambda$ as above and set $\newdelta := \newdelta(C,\newlambda_{0},\Delta,\Lambda^{c})$. Consider two copies $\newsigma_{\Lambda}^{\newomega}$ and $\newsigma_{\Lambda}^{\newtau}$ of $\newsigma_{\Lambda}$, which are distributed according to the measures $\newmu_{\Lambda}^{\newomega}$ and $\newmu_{\Lambda}^{\newtau}$ respectively. We claim that there is a coupling $\neweta$ of $\newsigma_{\Lambda}^{\newomega}$ and $\newsigma_{\Lambda}^{\newtau}$ such that for every $u_{\Delta} \varin E^{\Delta}$ with $\newmu(\newsigma_{\Delta} = u_{\Delta}) > 0$ we have
    \begin{align} \label{eq9}
        \neweta\big(\newsigma_{\Lambda}^{\newomega} \neq u_{\Delta}^{~}\!~|\!~\newsigma_{\Lambda}^{\newtau} = u_{\Delta}^{~}\big) \leq \newdelta/2
    \end{align}
    as well as
    \begin{align} \label{eq10}
        \neweta\big(\newsigma_{\Lambda}^{\newtau} \neq u_{\Delta}^{~}\!~|\!~\newsigma_{\Lambda}^{\newomega} = u_{\Delta}^{~}\big) \leq \newdelta/2.
    \end{align}
    Assuming that such a coupling $\neweta$ exists and observing that 
    \begin{align*}
        \frac{\newmu_{\Lambda}^{\newomega}(\newsigma_{\Delta} = u_{\Delta})}{\newmu_{\Lambda}^{\newtau}(\newsigma_{\Delta} = u_{\Delta})} &= \frac{\neweta(\newsigma_{\Lambda}^{\newomega} = u_{\Delta}^{~})}{\neweta(\newsigma_{\Lambda}^{\newtau} = u_{\Delta}^{~})} \\[0.8em] &= \frac{\neweta(\newsigma_{\Lambda}^{\newomega} = u_{\Delta}^{~} = \newsigma_{\Lambda}^{\newtau})}{\neweta(\newsigma_{\Lambda}^{\newtau} = u_{\Delta}^{~})}~\frac{\neweta(\newsigma_{\Lambda}^{\newomega} = u_{\Delta}^{~})}{\neweta(\newsigma_{\Lambda}^{\newomega} = u_{\Delta}^{~} = \newsigma_{\Lambda}^{\newtau})} \\[0.8em] &= \frac{\neweta(\newsigma_{\Lambda}^{\newomega} = u_{\Delta}^{~}\!~|\!~\newsigma_{\Lambda}^{\newtau} = u_{\Delta}^{~})}{  \neweta(\newsigma_{\Lambda}^{\newtau} = u_{\Delta}^{~}\!~|\!~\newsigma_{\Lambda}^{\newomega} = u_{\Delta}^{~})},
    \end{align*}
  we obtain
    \[1-\newdelta \leq \Big(1 - \frac{\newdelta}{2}\Big) \leq \frac{\newmu_{\Lambda}^{\newomega}(\newsigma_{\Delta} = u_{\Delta})}{\newmu_{\Lambda}^{\newtau}(\newsigma_{\Delta} = u_{\Delta})} \leq \Big(1 - \frac{\newdelta}{2}\Big)^{-1} \leq 1 + \newdelta\]
    by the fact that $\newdelta < 1$ and thus
    \[\bigg|\frac{\newmu_{\Lambda}^{\newomega}(\newsigma_{\Delta} = u_{\Delta})}{\newmu_{\Lambda}^{\newtau}(\newsigma_{\Delta} = u_{\Delta})} - 1\bigg| \leq \newdelta\]
    for all $u_{\Delta} \varin E^{\Delta}$, which gives (\ref{eq5}).

Accordingly it remains to construct the above coupling. To this end we decompose the set $\Lambda\setminus \Delta$ into the three sets $\Lambda_{1},\Lambda_{2},\Lambda_{3}$ defined by
        \begin{align*}
            \Lambda_{1} :&= \big\{s \in \Lambda\colon 0 < d(s,\Delta)/d(s,\Lambda^{c}) \leq 1/2\big\}, \\
            \Lambda_{2} :&= \big\{s \in \Lambda\colon 1/2 < d(s,\Delta)/d(s,\Lambda^{c}) \leq 2\big\}, \\
            &\Lambda_{3} := \big\{s \in \Lambda\colon d(s,\Delta)/d(s,\Lambda^{c}) > 2\big\}.
        \end{align*}
We start by showing that \begin{align} \label{eq6}
        \newpartial_{\!k}(\Delta \cup \Lambda_{1}) \subseteq \Lambda_{2},
        \end{align}
which will be used below. To verify (\ref{eq6}) note that by the choice of the sets $\Delta,\Lambda$ we have $\exp(\newlambda_{0} d(\Delta,\Lambda^{c})) > C$ and thus
        \begin{align*}
        d(\Delta,\Lambda^{c}) > \ln C /\newlambda_{0} \geq 3k\newlambda_{0}/\newlambda_{0} = 3k
        \end{align*}
by the choice of $C$. Now fix $t \varin \newpartial_{\!k}(\Delta \cup \Lambda_{1})$. Then $t \varin \Lambda_{1}^{c}$ and thus $d(t,\Delta)/d(t,\Lambda^{c}) > 1/2$. We distinguish two cases. If $d(t,\Delta) \leq 2k$, then the triangle inequality and the above inequality imply $d(t,\Lambda^{c}) \geq k$, so we obtain $d(t,\Delta)/d(t,\Lambda^{c}) \leq 2$ and thus $t \in \Lambda_{2}$. If on the other hand $d(t,\Delta) > 2k$ we may fix $s \in \Lambda_{1}$ with $d(s,t) \leq k$. Then again the triangle inequality implies $d(s,\Delta) \geq k$. Since $s \varin \Lambda_{1}$, this gives in turn $d(s,\Lambda^{c}) \geq 2d(s,\Delta) \geq 2k$. Using the triangle inequality again we obtain $d(t,\Delta) \leq d(s,\Delta) + k$ as well as $d(s,\Lambda^{c}) \leq d(t,\Lambda^{c}) + k$, so we get
    \begin{align*}
        d(t,\Lambda^{c}) \geq d(s,\Lambda^{c}) - k = d(s,\Lambda^{c}) - kd(s,\Lambda^{c})/d(s,\Lambda^{c}) \geq d(s,\Lambda^{c})/2
    \end{align*}
and
    \begin{align*}
        d(t,\Delta) \leq d(s,\Delta) + k = d(s,\Delta) + k d(s,\Delta)/d(s,\Delta) \leq 2 d(s,\Delta).
    \end{align*}
Both together gives
    \[\frac{d(t,\Delta)}{d(t,\Lambda^{c})} \leq \frac{4}{2} = 2\]
by the fact that $s \varin \Lambda_{1}$, which implies $t \varin \Lambda_{2}$.

We will now describe the coupling $\neweta$. By the choice of $\newomega,\newtau$ we have
    \begin{align}
        \big\|\smash{\newsigma_{\Lambda_{2}}}_{*}\newmu_{\Lambda}^{\newomega} - \smash{\newsigma_{\Lambda_{2}}}_{*}\newmu_{\Lambda}^{\newtau}\|^{~}_{\text{Var}} \leq \newdelta(C_{1},\newlambda,\Lambda_{2},\Lambda^{c}) := \newdelta_{1}.
    \end{align}
By the maximal coupling lemma this implies that there exists a coupling $\neweta_{0}$ of $\smash{\newsigma_{\Lambda_{2}}}_{*}\newmu_{\Lambda}^{\newomega}$ and $\smash{\newsigma_{\Lambda_{2}}}_{*}\newmu_{\Lambda}^{\newtau}$ such that
    \[\neweta_{0}\big(\newsigma^{\omega}_{\Lambda_{2}} \neq \newsigma^{\tau}_{\Lambda_{2}}\big) \leq \newdelta_{1}.\]
Based on this we choose a pair $(v_{\Lambda},w_{\Lambda})$ of configurations as follows. In the first step we choose a pair $(v_{\Lambda_{2}},w_{\Lambda_{2}})$ according to $\neweta_{0}$. In the second step we choose $v_{\Lambda_{3}}$ and $w_{\Lambda_{3}}$ independently under $\smash{\newsigma_{\Lambda_{3}}}_{*}\newmu_{\Lambda}^{\newomega}(\!\!~\cdot~\!\!|\newsigma_{\Lambda_{2}} = v_{\Lambda_{2}})$ and $\smash{\newsigma_{\Lambda_{3}}}_{*}\newmu_{\Lambda}^{\tau}(\!\!~\cdot~\!\!|\newsigma_{\Lambda_{2}} = w_{\Lambda_{2}})$ respectively. Now, if $v_{\Lambda_{2}} = w_{\Lambda_{2}}$, we choose $v_{\Delta \cup \Lambda_{1}}$ under $\smash{\newsigma_{\Delta \cup \Lambda_{1}}}_{*}\newmu_{\Lambda}^{\newomega}(\!\!~\cdot~\!\!|\newsigma_{\Lambda_{2} \cup \Lambda_{3}} = v_{\Lambda_{2} \cup \Lambda_{3}})$ and set $w_{\Delta \cup \Lambda_{1}} = v_{\Delta \cup \Lambda_{1}}$. Otherwise, we choose $v_{\Delta \cup \Lambda_{1}}$ and $w_{\Delta \cup \Lambda_{1}}$ independently under $\smash{\newsigma_{\Delta \cup \Lambda_{1}}}_{*}\newmu_{\Lambda}^{\newomega}(\!\!~\cdot~\!\!|\newsigma_{\Lambda_{2} \cup \Lambda_{3}} = v_{\Lambda_{2} \cup \Lambda_{3}})$ and $\smash{\newsigma_{\Delta \cup \Lambda_{1}}}_{*}\newmu_{\Lambda}^{\tau}(\!\!~\cdot~\!\!|\newsigma_{\Lambda_{2} \cup \Lambda_{3}} = v_{\Lambda_{2} \cup \Lambda_{3}})$ respectively. 

We may check that $\neweta$ defines indeed a coupling as follows. Since in the first step we use a coupling of $\smash{\newsigma_{\Lambda_{2}}}_{*}\newmu_{\Lambda}^{\newomega}$ and $\smash{\newsigma_{\Lambda_{2}}}_{*}\newmu_{\Lambda}^{\newtau}$, we obtain in the second step a coupling of $\smash{(\newsigma_{\Lambda_{2} \cup \Lambda_{3}})}_{*}\newmu_{\Lambda}^{\newomega}$ and $\smash{(\newsigma_{\Lambda_{2} \cup \Lambda_{3}})}_{*}\newmu_{\Lambda}^{\newtau}$. Now note that for any two configurations $v_{\Lambda_{2} \cup \Lambda_{3}}$ and $w_{\Lambda_{2} \cup \Lambda_{3}}$ with $v_{\Lambda_{2}} = w_{\Lambda_{2}}$ we have 
    \[\smash{(\newsigma_{\Delta \cup \Lambda_{1}})}_{*}\newmu_{\Lambda}^{\newomega}(\!\!~\cdot~\!\!|\newsigma_{\Lambda_{2} \cup \Lambda_{3}} = v_{\Lambda_{2} \cup \Lambda_{3}}) = \smash{(\newsigma_{\Delta \cup \Lambda_{1}})}_{*}\newmu_{\Lambda}^{\tau}(\!\!~\cdot~\!\!|\newsigma_{\Lambda_{2} \cup \Lambda_{3}} = v_{\Lambda_{2} \cup \Lambda_{3}})\]
by \ref{eq6}. Using this observation it is easy to see that we obtain in the third step a coupling of $\newmu_{\Lambda}^{\newomega}$ and $\newmu_{\Lambda}^{\newtau}$ as desired.

To see that $\neweta$ satisfies the claimed properties let us fix a configuration $u_{\Delta} \in E^{\Delta}$. Then we may compute
    \begin{align*}
    \neweta\big(\newsigma^{\newtau}_{\Delta} \neq u^{~}_{\Delta}~\!\!\big|\!\!~\newsigma^{\newomega}_{\Delta} = &u^{~}_{\Delta}\big) = 
\neweta(\newsigma_{\Lambda_{2}}^{\newomega} = \newsigma_{\Lambda_{2}}^{\newtau}~\!\!|\!\!~\newsigma^{\newomega}_{\Delta} = u^{~}_{\Delta})~      \neweta(\newsigma^{\newtau}_{\Delta} \neq u^{~}_{\Delta}~\!\!\big|\!\!~\newsigma^{\newomega}_{\Delta} = u^{~}_{\Delta},\newsigma_{\Lambda_{2}}^{\newomega} = \newsigma_{\Lambda_{2}}^{\newtau})~+ \vphantom{\int} \\
    &+~\neweta(\newsigma_{\Lambda_{2}}^{\newomega} \neq \newsigma_{\Lambda_{2}}^{\newtau}~\!\!|\!\!~\newsigma^{\newomega}_{\Delta} = u^{~}_{\Delta})~      \neweta(\newsigma^{\newtau}_{\Delta} \neq u^{~}_{\Delta}~\!\!\big|\!\!~\newsigma^{\newomega}_{\Delta} = u^{~}_{\Delta},\newsigma_{\Lambda_{2}}^{\newomega} \neq \newsigma_{\Lambda_{2}}^{\newtau}).
    \end{align*}
By the definition of the coupling $\neweta$ we also have
    \[\neweta\big(\newsigma^{\newtau}_{\Delta} \neq u^{~}_{\Delta}~\!\!\big|\!\!~\newsigma^{\newomega}_{\Delta} = u^{~}_{\Delta},\newsigma^{\newomega}_{\Lambda_{2}} = \newsigma^{\newtau}_{\Lambda_{2}}\big) = 0,\]
so we obtain
    \[\neweta\big(\newsigma^{\newtau}_{\Delta} \neq u^{~}_{\Delta}~\!\!\big|\!\!~\newsigma^{\newomega}_{\Delta}= u^{~}_{\Delta}\big) \leq  \neweta(\newsigma_{\Lambda_{2}}^{\newomega} \neq \newsigma_{\Lambda_{2}}^{\newtau}~\!\!|\!\!~\newsigma^{\newomega}_{\Delta} = u^{~}_{\Delta}).\]
Accordingly it remains to estimate the latter expression. To this end we consider the random variable $g\colon \Omega \times \Omega \to \mathbb{R}$ defined by 
    \[g(\newxi,\zeta) := \neweta(\newsigma^{\newomega}_{\Lambda_{2}} \neq \newsigma^{\newtau}_{\Lambda_{2}}|\newsigma^{\newomega}_{\Lambda_{2}} = \newxi^{~}_{\Lambda_{2}})\]
for $\newxi \varin \Omega$ and set $c := d(\Lambda_{2},\Lambda^{c})/3$. By Lemma \ref{lemma1} we have
    \[\big|\newmu_{\Lambda}^{\newomega}\big(g > e^{-\newlambda c}|\newsigma_{\Delta} = u_{\Delta}\big) - \newmu_{\Lambda}^{\newomega}\big(g > e^{-\newlambda c})\big| \leq \newdelta(C_{3},\newlambda,\Lambda_{2},\Delta \cup \Lambda^{c}\big) =: \newdelta_{2}.\]
Moreover obtain
    \[\mathbb{E}[g] = \neweta(\newsigma^{\newomega}_{\Lambda_{2}} \neq \newsigma^{\newtau}_{\Lambda_{2}}) =\neweta_{0}(\newsigma^{\newomega}_{\Lambda_{2}} \neq \newsigma^{\newtau}_{\Lambda_{2}}),\]
which, by Markov's inequality, implies
    \begin{align*}
        \newmu_{\Lambda}^{\newomega}(g > e^{-\newlambda c}) &= \neweta(g > e^{-\newlambda c}) \leq e^{\newlambda c} \mathbb{E}[g]  \vphantom{\sum_{~}} 
 \leq e^{\newlambda c} \newdelta_{1} \\ &= C_{1}\sum_{t \in \Lambda_{2}}\sum_{s \in \Lambda^{c}}e^{\newlambda(c - d(s,t))} \\
&\leq C_{1} \sum_{t \in \Lambda_{2}}\sum_{s \in \Lambda^{c}} e^{-2\newlambda d(s,t)/3} \\ &= \newdelta(C_{1},2\lambda/3,\Lambda_{2},\Lambda^{c}) =: \newdelta_{3}.
    \end{align*}
Both together gives
    \begin{align} \label{eq7}
        \nonumber \neweta(g > e^{-\newlambda c}|\newsigma^{\newomega}_{\Delta} = u^{~}_{\Delta}) &= \newmu_{\Lambda}^{\newomega}(g > e^{-\newlambda c}|\newsigma_{\Delta} = u_{\Delta}) \\ &\leq \newmu_{\Lambda}^{\newomega}(g > e^{-\newlambda c}) + \newdelta_{2} \leq \newdelta_{3} + \newdelta_{2}.
    \end{align}
On the other hand we get
    \begin{align*}
        \neweta(\newsigma^{\newomega}_{\Lambda_{2}} \neq \newsigma^{\newtau}_{\Lambda_{2}},~g \leq e^{-\newlambda c}&|\newsigma^{\newomega}_{\Delta} = u^{~}_{\Delta}) \leq \neweta(\newsigma^{\newomega}_{\Lambda_{2}} \neq \newsigma^{\newtau}_{\Lambda_{2}}~\!\!|~\! g \leq e^{-\newlambda  c},\newsigma^{\newomega}_{\Delta} = u^{~}_{\Delta}) \\[0.8em]
&\leq \max\big\{\neweta(\newsigma^{\newomega}_{\Lambda_{2}} \neq \newsigma^{\newtau}_{\Lambda_{2}}|  \newsigma^{\newomega}_{\Lambda_{2}} = \newxi_{\Lambda_{2}}, \newsigma^{\newomega}_{\Delta} = u_{\Delta})\colon g(\newxi,\zeta) \leq e^{-\newlambda c}\big\}.
    \end{align*}
It is not hard to check that by the construction of the coupling $\neweta$ the random variable $\newsigma^{\newomega}_{\Delta}$ is conditionally independent of the random variable $\newsigma^{\newtau}_{\Lambda_{2}}$ given $\newsigma^{\newomega}_{\Lambda_{2}}$. But conditional independence is a symmetric relation, so we obtain as well that $\newsigma^{\newtau}_{\Lambda_{2}}$ is conditionally independent of $\newsigma^{\newomega}_{\Delta}$ given $\newsigma^{\newomega}_{\Lambda_{2}}$. Thus we get
    \[\neweta(\newsigma^{\newomega}_{\Lambda_{2}} \neq  \newsigma^{\newtau}_{\Lambda_{2}}|  \newsigma^{\newomega}_{\Lambda_{2}} = \newxi_{\Lambda_{2}}, \newsigma^{\newomega}_{\Delta} = u^{~}_{\Delta}) = \neweta(\newsigma^{\newomega}_{\Lambda_{2}} \neq  \newsigma^{\newtau}_{\Lambda_{2}}|  \newsigma^{\newomega}_{\Lambda_{2}} = \newxi_{\Lambda_{2}})\]
for all $\newxi \in \Omega$, which by the definition of $g$ implies
    \begin{align} \label{eq8}
        \nonumber \neweta(\newsigma^{\newomega}_{\Lambda_{2}} \neq \newsigma^{\newtau}_{\Lambda_{2}}, ~&g \leq e^{-\newlambda c}|\newsigma^{\newomega}_{\Delta} = u^{~}_{\Delta}) \leq e^{-\newlambda c} = e^{-\newlambda d(\Lambda_{2},\Lambda^{c})/3} \vphantom{\sum_{~}} \\
&\leq \sum_{s \in \Lambda_{2}}\sum_{t \in \Lambda^{c}} e^{-\newlambda d(s,t)/3} = \newdelta(1,\newlambda/3,\Lambda_{2},\Lambda^{c}) =: \newdelta_{4}.
    \end{align} 
The equalities (\ref{eq7}) and (\ref{eq8}) together imply
    \begin{align*}
        \neweta(\newsigma^{\newomega}_{\Lambda_{2}} \neq \newsigma^{\newtau}_{\Lambda_{2}}|\newsigma^{\newomega}_{\Delta} = u^{~}_{\Delta}) &= \neweta(\newsigma^{\newomega}_{\Lambda_{2}} \neq \newsigma^{\newtau}_{\Lambda_{2}}, g \leq e^{-\newlambda c}| \newsigma^{\newomega}_{\Delta} = u^{~}_{\Delta}) \\ &~~~+ \neweta(\newsigma^{\newomega}_{\Lambda_{2}} \neq \newsigma^{\newtau}_{\Lambda_{2}}, g > e^{-\newlambda c}|\newsigma^{\newomega}_{\Delta} = u^{~}_{\Delta}) \\
		&\leq \neweta(\newsigma^{\newomega}_{\Lambda_{2}} \neq \newsigma^{\newtau}_{\Lambda_{2}}, g \leq e^{-\newlambda c}|\newsigma^{\newomega}_{\Delta} = u^{~}_{\Delta})
		\\ &~~~+ \neweta(g > e^{-\newlambda c}|\newsigma^{\newomega}_{\Delta} = u^{~}_{\Delta}) \leq \newdelta_{2} + \newdelta_{3} + \newdelta_{4}.
    \end{align*}
Now fix $s \varin \Delta$ and $t \in \Lambda_{2}$. Then by the triangle inequality we obtain
    \[d(s,t) \geq d(t,\Delta) \geq d(t,\Lambda^{c})/2 \geq (d(s,\Lambda^{c}) - d(s,t))/2,\]
which implies $d(s,t) \geq d(s,\Lambda^{c})/3$. Similarly we can show that for $s \varin \Lambda^{c}$ and $t \in \Lambda_{2}$ we have $d(s,t) \geq d(s,\Delta)/3$. Thus we have
    \begin{align*}
        \Sigma(\newlambda/3,&\Lambda_{2},\Delta \cup \Lambda^{c}) = \sum_{t \!\!\varin\!\! \Lambda_{2}}\sum_{s \!\!\varin\!\! \Delta \cup \Lambda^{c}} e^{-\newlambda d(s,t)/3} \\
&\leq \sum_{s \!\!\varin\!\! \Delta}~\sum_{d(s,t) \geq d(s,\Lambda^{c})/3} e^{-\newlambda d(s,t)/3} + \sum_{s \!\!\varin\!\! \Lambda^{c}}~\sum_{d(s,t) \geq d(s,\Delta)/3} e^{-\newlambda d(s,t)/3} \\
&= \sum_{s \!\!\varin\!\! \Delta}~\sum_{r \geq d(s,\Lambda^{c})/3} |B_{r}(s)| e^{-\newlambda r/3} + \sum_{s \!\!\varin\!\! \Lambda^{c}}~\sum_{r \geq d(s,\Delta)/3} |B_{r}(s)| e^{-\newlambda r/3} \\
&\leq  \sum_{s \!\!\varin\!\! \Delta}~\sum_{r \geq d(s,\Lambda^{c})/3} e^{r(\newtheta - \newlambda/3)} +  \sum_{s \!\!\varin\!\! \Lambda^{c}}~\sum_{r \geq d(s,\Delta)/3} e^{r(\newtheta - \newlambda/3)} \\
&\leq   \newbeta \sum_{s \in \Delta}e^{(\newtheta - \newlambda/3)d(s,\Lambda^{c})/3} +   \newbeta \sum_{s \in \Lambda^{c}} e^{(\newtheta - \newlambda/3) d(s,\Delta)/3} \\
&\leq  \newbeta \sum_{s \in \Delta} \sum_{t \in \Lambda^{c}}e^{-\newlambda_{0}d(s,t)} + \newbeta \sum_{s \in \Lambda^{c}}\sum_{t \in \Delta} e^{-\newlambda_{0}d(s,t)} \\
&= 2\newbeta \Sigma(\newlambda_{0},\Delta, \Lambda^{c}).
    \end{align*}
Noting that $\newdelta_{j} \leq  \newdelta(C_{4}, \newlambda/3,\Lambda_{2},\Delta \cup \Lambda^{c})$ for $j = 2,3,4$ this shows that
    \[\newdelta_{2} + \newdelta_{3} + \newdelta_{4} \leq 6 C_{4} \newbeta \Sigma(\newlambda_{0},\Delta,\Lambda^{c}) \leq \newdelta(C,\newlambda_{0},\Delta,\Lambda^{c})/2 = \newdelta/2, \]
which proves $(\ref{eq9})$. To show $(\ref{eq10})$ we consider the random variable $h\colon \Omega \times \Omega \to \mathbb{R}$ given by 
    \[h(\newxi,\newzeta) := \neweta(\newsigma^{\newomega}_{\Lambda_{2}} \neq \newsigma^{\newtau}_{\Lambda_{2}}| \newsigma^{\newtau}_{\Delta} = \newzeta^{~}_{\Delta})\]
for $(\newxi,\newzeta) \varin \Omega \times \Omega$ and note that by (\ref{eq6}), the fact that $\newmu$ is $k$-Markovian and the construction of the coupling $\neweta$ the random variables $\newsigma^{\newtau}_{\Delta}$ and $\newsigma^{\newomega}_{\Lambda_{2}}$ are conditionally independent given $\newsigma^{\tau}_{\Lambda_{2}}$. This observation allows us to argue as before to verify the claim. \hfill $\Box$

\section{Applications} \label{applications}

Typical examples of exponential mixing arise in the theory of Gibbs fields. As we shall recall below it is well known that for Gibbs fields over general graphs a sufficiently small Dobrushin constant implies exponential $\newupphi$-mixing with arbitrarily high rate of decay. Based on our results this fact will enable us to formulate a general criterion for exponential $\newuppsi$-mixing of such Gibbs fields in terms of their Dobrushin constant. This allows us to obtain examples of exponentially $\newuppsi$-mixing Gibbs fields over more general graphs than square lattices. We shall illustrate this in the case of the Ising and Potts model on a regular tree.

We start by recalling some basic facts about Gibbs fields. For a detailed exposition of the topic we refer to \cite{Geo11}. A \textit{Gibbs field} is a random field that realizes a certain specification of conditional probabilities. Formally such a specification is given by a family $(\newrho_{\Lambda})_{\Lambda \sqsubseteq S}$ of probability kernels $\newrho_{\Lambda}$ from $\mathscr{F}_{\Lambda^{c}}$ to $\mathscr{F}_{\Lambda}$ satisfying the consistency condition
	\[\newrho_{\Lambda}^{\newomega}(A) = \sum_{w_{\Lambda} \!\!\varin\!\!E^{\Lambda}}\newrho_{\Delta}^{w_{\Lambda}\newomega_{\Lambda^{c}}}(A)\]
for all $\Delta \subseteq \Lambda \sqsubseteq S$, $\newomega \varin \Omega$ and $A \varin \mathscr{F}_{\Delta}$. By a \textit{realization} of $\newrho$ we mean a probability measure $\newmu$ such that for all $\Lambda \sqsubseteq S$ and all $A \varin \mathscr{F}_{\Lambda}$ we have
	\[\newmu_{\Lambda}^{\newomega}(A) = \newrho_{\Lambda}^{\newomega}(A)\]
for $\newmu$-almost every $\newomega \in \Omega$. Similarly as in the the case of probability measures we will call a specification $\newrho$ \textit{$k$-Markovian} for a given $k \in \mathbb{N}_{0}$ if $\newrho^{\newomega}_{\newpartial_{k}\Lambda} = \newrho^{\newtau}_{\newpartial_{k}\Lambda}$ implies $\newrho_{\Lambda}^{\newomega} = \newrho_{\Lambda}^{\tau}$ for all $\Lambda \sqsubseteq S$ and all $\newomega,\newtau \varin \Omega$. Obviously any realization of a $k$-Markovian specification is $k$-Markovian.  

Gibbs measures appear as realizations of specifications given by potentials. A \textit{potential} is given by a family $\Phi  = (\Phi_{\Delta})_{\Delta \sqsubseteq S}$ of $\mathscr{F}_{\Delta}$-measurable functions $\Phi_{\Delta}\colon \Omega \to \mathbb{R}$ such that the series
    \[H_{\Lambda}^{\Phi}(\newomega) := \sum_{\Delta \cap \Lambda \neq \emptyset} \Phi_{\Delta}(\newomega)\]
converges absolutely for every $\Lambda\sqsubseteq S$ and all $\newomega \varin \Omega$. The arising map $H_{\Lambda}^{\Phi}\colon \Omega \to \mathbb{R}$ is called the \textit{Hamiltonian} of $\Phi$ at $\Lambda$, while the \textit{partition function} $Z^{\Phi}_{\Lambda}\colon\Omega \to \mathbb{R}$ of $\Phi$ at $\Lambda$ is defined by
    \[Z_{\Lambda}^{\Phi}(\newomega) := \sum_{w_{\Lambda} \!\!\varin\!\! E^{\Lambda}} e^{-H^{\Phi}_{\Lambda}(w_{\Lambda}\newomega_{\Lambda^{c}})}\]
for $\omega \varin \Omega$. Given a potential $\Phi$ there is a unique specification $^{\Phi\!\!}\newrho$ satisfying
    \[(^{{\Phi}\!\!}\newrho)_{\Lambda}^{\newomega}\big(\newsigma_{\Lambda} = w_{\Lambda}\big) = \frac{e^{-H_{\Lambda}^{\Phi}(w_{\Lambda}\newomega_{\Lambda^{c}})}}{Z^{\Phi}_{\Lambda}(\newomega)}\]
for every $\Lambda \sqsubseteq S$, $\newomega \varin \Omega$ and $w_{\Lambda} \varin E^{\Lambda}$. We refer to a realization of $^{\Phi\!\!}\newrho$ as a \textit{Gibbs measure} with potential $\Phi$. The above formula immediately implies that Gibbs measures are non-degenerate. 

A common assumption on potentials considered in examples is that they have finite range, where the \textit{range} of a potential $\Phi$ is defined by
\[r(\Phi) := \sup\big\{\text{diam}(\Lambda)\colon \Phi_{\Lambda} \neq 0\big\}.\]
Finite range potentials correspond to Markovian specifications. In fact, it is not difficult to see that a specification and thus every Gibbs measure arising from a potential of range $k$ is $k$-Markovian.

Gibbs measures realizing a given potential $\Phi$ exist under very general conditions. A mild sufficient condition is \textit{absolute summability} of the potential $\Phi$, which requires that for all $s \varin S$ we have
    \[\sum_{\Lambda \!\varni\! s} \|\Phi_{\Lambda}\|_{\infinity} < \infinity.\]
Obviously finite range potentials are absolutely summable, so they admit at least one realizing Gibbs measure. A key problem in the theory of Gibbs fields is the study of phase transition, i.e. the question under which conditions there is more than one realizing Gibbs measure. A well known condition ensuring uniqueness of Gibbs measures is due to Dobrushin, see \cite{Dob68}. Given a specification $\newrho$ let us consider the dependency matrix $\Gamma = (\newgamma_{st})_{s,t \!\varin\! S}$ defined by
    \[\newgamma_{st} := \sup\Big\{\|\newrho_{\{s\}}^{\newomega} - \newrho_{\{s\}}^{\newtau}\|_{\text{Var}}\colon \newomega_{\{t\}^{c}} = \newtau_{\{t\}^{c}} \Big\}\]
for $s,t \varin S$. The so called \textit{Dobrushin constant} $\newgamma$ of $\Phi$ is then given by
    \[\newgamma := \sup_{s \!\varin\! S}\sum_{t \!\varin\! S} \newgamma_{st}.\]
Dobrushin's criterion states that if $\newgamma < 1$ there is at most one Gibbs measure realizing $\newrho$. For many classes of potentials one knows explicit estimates of Dobrushin's constant some of which will be discussed below. Moreover, it is well known that a small Dobrushin constant ensures $\newupphi$-mixing with exponential rate of decay. For the sake of completeness let us recall this result, which goes also back to Dobrushin, see \cite{Dob68} and \cite{Dob70}.

\begin{theorem} \thm \label{criterion_phi-mixing}
    Fix $k \in \mathbb{N}_{0}$ and let $\newrho$ be a $k$-Markovian specification with Dobrushin constant $\newgamma$. If $\newgamma < \exp(-3\newtheta k)$, then the unique Gibbs measure $\newmu$ realizing $\newrho$ is $\newupphi$-mixing with exponential decay rate $\newlambda(k,\newgamma) := - \ln(\newgamma)/k$.
\end{theorem}

\textit{Proof:} Since $\newrho$ is $k$-Markovian the entries $\newgamma_{st}$ of the dependency matrix $\Gamma$ vanish whenever $d(s,t) > k$. Using this it is not difficult to see that for every $n \in \mathbb{N}$ the $n$-th power $\Gamma^{n}$ of $\Gamma$, whose entries shall be denoted by $\smash{\smash{\newgamma_{st}^{(n)}}}$, is well-defined. Assume that $\newgamma < 1$ is satisfied and let $\newmu$ be the unique realization of $\newrho$. Then it is easy to check that we have $\smash{\newgamma_{st}^{(n)}} \leq \newgamma^{\!~n}$ for all $n \varin \mathbb{N}$ and $s,t \varin S$. Thus, setting $\Gamma^{0} = \mathrm{Id}$, we may consider the matrix $Q := \sum_{n=0}^{\infinity} \Gamma^{n}$, whose entries we shall denote by $q_{st}$. It is well known that for all $\Delta \subseteq \Lambda \sqsubseteq S$ we have the upper bound
    \begin{align} \label{eq1}
        \sup\Big\{|\newrho_{\Lambda}^{\newomega}(A) - \newmu(A)|\colon A \in \mathscr{F}_{\Delta}\Big\} \leq \sum_{s \in \Delta, t \in \Lambda^{c}}q_{st}
    \end{align}
for all $\newomega \varin \Omega$, see \cite[Theorem 8.23]{Geo11}. Using once again that $\newrho$ is $k$-Markovian it is easy to check that $\smash{\newgamma_{st}^{(n)}} = 0$ whenever $d(s,t) > nk$. Accordingly we get
    \begin{align*}
        q_{st} = \sum_{n \geq d(s,t)/k}\newgamma_{st}^{(n)} \leq \sum_{n \geq d(s,t)/k} \newgamma^{~\!n} \leq \frac{1}{1{-}\newgamma}~ \newgamma^{~\!d(s,t)/k} = \frac{1}{1{-}\newgamma}~e^{d(s,t)\ln (\newgamma)  /k}.
    \end{align*}
Combining this with (\ref{eq1}) gives
    \[|\newmu^{\newomega}_{\Lambda}(A) - \newmu^{\newtau}_{\Lambda}(A)| \leq |\newrho_{\Lambda}^{\newomega}(A) - \newmu(A)| + |\newmu(A) - \newrho_{\Lambda}^{\newtau}(A)| \leq C~\Sigma(\newlambda,\Delta,\Lambda^{c})\]
with $C = 2/(1{-}\newgamma)$ and $\newlambda(k,\newgamma) = - \ln(\newgamma)/k$ for all $\Delta \subseteq \Lambda \sqsubseteq S$, all $A \varin \mathscr{F}_{\Delta}$ and all $\newomega,\newtau \varin \Omega$. Thus, by Proposition \ref{prop1}, we obtain that $\newmu$ is exponential $\newupphi$-mixing with decay rate $\newlambda(k,\newgamma)$ whenever $\newgamma < \exp(-3\newtheta k)$. \hfill $\Box$ \\

We can use the above theorem in combination with Theorem \ref{main} to obtain the following criterion for exponential $\newuppsi$-mixing of $k$-Markovian Gibbs fields in terms of the corresponding Dobrushin constant.

\begin{theorem} \corollary \label{criterion_psi-mixing}
    Fix $k \in \mathbb{N}_{0}$ and let $\newrho$ be a $k$-Markovian specification with Dobrushin constant $\newgamma$. If $\newgamma < \exp(-3\newtheta k)$, then the unique Gibbs measure $\newmu$ realizing $\newrho$ is $\newuppsi$-mixing with exponential decay rate 
        \[\newlambda(k,\newgamma,\newtheta) := -\frac{1}{9} \left(\frac{\ln\newgamma}{k} + 3\newtheta\right).\]
\end{theorem}

\medskip

In the remainder of this section we will apply the above criterion to Gibbs fields on regular trees. To this end fix $d > 2$ and let $\mathscr{T}_{d}$ denote the $d$-regular tree. Note that the exponential growth rate of $\mathscr{T}_{d}$ is bounded by $\ln d$. We shall study instances of exponential $\newuppsi$-mixing in two classical regimes of uniqueness: The Ising model at high temperatures and the antiferromagnetic Potts model with large number of states. 

To this end let us briefly recall how these models are defined in the present context. Fix $\newbeta > 0$ and set $E := \{-1,1\}$. Then the \textit{Ising model} at \textit{inverse temperature} $\newbeta$ corresponds to the potential $\Phi^{\newbeta}$ defined by $\smash{\Phi^{\newbeta}_{\Lambda}}(\newomega) := -\newbeta\newomega_{s}\newomega_{t}$ for sets of the form $\Lambda = \{s,t\}$ with $s \sim t$ and $\smash{\Phi^{J,h}_{\Lambda}} :\equiv 0$ for all other sets $\Lambda$. Thus the respective Hamiltonian $H^{\newbeta}$ takes the form
        \[H^{\newbeta}_{\Lambda}(\newomega) := -\newbeta\sum_{t \!\varin\! \Lambda}\sum_{s \sim t}\newomega_{s}\newomega_{t}.\]
    Now set $E := \{1,...,N\}$. Then the \textit{antiferromagnetic $N$-state Potts model} at \textit{inverse temperature} $\newbeta$ corresponds to the potential $\Psi^{\newbeta}$ defined by $\Psi^{\newbeta}_{\Lambda}(\newomega) := \newbeta\newdelta(\newomega_{s},\newomega_{t})$ for $\Lambda = \{s,t\}$ with $s \sim t$, where $\newdelta$ denotes the Kronecker delta, and $\Psi^{J}_{\Lambda} :\equiv 0$ for every other set $\Lambda$.
    
    The first corollary shows that below a certain inverse temperature $\newbeta$ the Ising model admits a unique Gibbs measure, which is exponentially $\newuppsi$-mixing. Moreover, making the inverse temperature sufficiently low we can make the rate of mixing arbitrarily high.

    \begin{theorem} \corollary
        Fix $\newbeta < \mathrm{artanh}(1/d^{4})$. Then the Ising model on $\mathscr{T}_{d}$ at inverse temperature $\newbeta$ admits a unique realizing Gibbs measure $\newmu$, which is $\newuppsi$-mixing with exponential decay rate given by
            \[-\frac{1}{9}\ln \tanh \newbeta - \frac{4}{9} \ln d.\]
    \end{theorem}

\textit{Proof:} In the above setting we may bound the Dobrushin constant $\newgamma$ from above by 
		\[\newgamma \leq d \tanh\newbeta,\]
cf. \cite[Example 8.9.2]{Geo11}. Thus Dobrushin's criterion implies uniqueness of the realizing Gibbs measure, whenever $\newbeta < \mathrm{artanh}(1/d)$. Since the potential $\Phi^{\newbeta}$ has range $1$, the unique realizing Gibbs measure $\newmu$ is a Markov measure.

Now let us assume that $\newbeta < \mathrm{artanh}(1/d^{4})$.  
Then we obtain $\newgamma < 1/d^{3} = \exp(-2\ln d)$. Since the exponential growth rate of $\mathscr{T}_{d}$ is smaller than $\ln d$ we may apply Corollary \ref{criterion_psi-mixing} to deduce that $\newmu$ is $\newuppsi$-mixing with exponential decay rate at least
	\[-\frac{1}{9}\ln \tanh \newbeta - \frac{4}{9} \ln d.\]
This finishes the proof. \hfill $\Box$ 

\medskip

The second corollary shows that above a certain number of states the antiferromagnetic Potts model admits a unique Gibbs measure, which is exponentially $\newuppsi$-mixing, independent of the inverse temperature. Moreover, the decay rate becomes arbitrarily high for sufficiently large number of states. 

\begin{theorem} \corollary
    Fix $N > 2d^{4} + d$. Then for any $\newbeta > 0$ the antiferromagnetic $N$-state Potts model on $\mathscr{T}_{d}$ admits a unique realizing Gibbs measure $\newmu$, which is $\newuppsi$-mixing with exponential decay rate given by 
        \[\frac{1}{9}\ln(N-d) - \frac{4}{9}\ln d - \frac{\ln 2}{9}.\]
\end{theorem}

\textit{Proof:} In the above setting we may estimate Dobrushin's constant from above by 
    \[\newgamma \leq \frac{2d}{N-d},\]
see \cite[Example 8.13.3]{Geo11}. Accordingly, for $N > 3d$, there is a unique Gibbs measure $\newmu$ corresponding to $\Psi^{\newbeta}$ for \textit{every} $\newbeta > 0$, which is Markov by the fact that $\Psi^{\newbeta}$ has range $1$. Furthermore, for $N > 2d^{4} + d$, we obtain $\newgamma < \exp(-3\ln d)$. Thus we may apply Corollary \ref{criterion_psi-mixing} to obtain that $\newmu$ is $\newuppsi$-mixing with exponential decay rate at least
	\[\frac{1}{9}\ln(N-d) - \frac{4}{9}\ln d - \frac{\ln 2}{9},\]
which shows the claim. \hfill $\Box$
\end{spacing}

\smallskip

\textsc{Mathematical Institute, University of Leipzig}

\textsc{Augustusplatz 10, 04109 Leipzig}

\textit{elias.zimmermann@math.uni-leipzig.de}

\end{document}